\documentclass[10pt]{article}
\usepackage[english]{babel}
\usepackage[utf8]{inputenc}
\usepackage[T1]{fontenc}
\usepackage{amsmath,amssymb,amsfonts,amsopn,amsthm}
\usepackage[a4paper, total={6in, 10in}]{geometry}

\usepackage{algorithm}
\usepackage{algpseudocode}

\usepackage[dvipsnames]{xcolor}
\usepackage{graphicx,epstopdf} 
\usepackage[labelformat=simple]{subcaption}
\usepackage[format=plain,
            textfont=it,
            size=footnotesize,labelsep=period]{caption}
\usepackage{pgfplots}
\usepackage{siunitx,array,booktabs}
\usepackage{ifpdf} 
\usetikzlibrary{patterns}
\pgfdeclarelayer{foreground layer} 
\pgfsetlayers{main,foreground layer}
\usepackage{commath}
\usepackage{bm}
\usepackage{enumitem}
\usepackage[colorinlistoftodos,italian]{todonotes}
\usepackage[colorlinks = true,linkcolor = black,
            citecolor = ForestGreen, anchorcolor = black]{hyperref}
\usepackage[capitalise]{cleveref}

\newcommand{\email}[1]{\href{#1}{#1}}

\pgfplotsset{compat=1.11}

\numberwithin{equation}{section}

\crefname{assumption}{Assumption}{Assumptions}
\crefname{remark}{Remark}{Remarks}
\crefname{example}{Example}{Examples}

\newcommand{\plotlinewidth}{1}

\newcommand{\aspectratio}{0.66}
\newcommand{\plotmarksizeu}{2.5}

\newcommand{\plotimscale}{0.98}
\newcommand{\plotimsizeu}{0.6}

\newcommand{\subfigsize}{0.49}

\newcommand{\legendfontscale}{1}
\newcommand{\legendmarkscale}{1}

\usepackage{calligra}
\usepackage{csvsimple}

\pgfplotsset{select coords between index/.style 2 args={
    x filter/.code={
        \ifnum\coordindex<#1\fi
        \ifnum\coordindex>#2\fi
    }
}}

\newcommand{\Nb}{\mathbb{N}}

\newcommand{\Rb}{\mathbb{R}}

\newcommand{\Rn}{\mathbb{R}^{n}}

\newcommand{\dt}{\partial_t}
\newcommand{\bigo}[1]{\mathcal{O}\left( #1 \right)}

\newcommand\blfootnote[1]{%
  \begingroup
  \renewcommand\thefootnote{}\footnote{#1}%
  \addtocounter{footnote}{-1}%
  \endgroup
}

\newcommand{\rhof}{\rho_F}
\newcommand{\rhos}{\rho_S}
\newcommand{\ff}{f_F}
\newcommand{\fs}{f_S}
\newcommand{\ARKC}{\text{ARKC}}

\title{Instabilities and order reduction phenomenon of an interpolation based multirate Runge--Kutta--Chebyshev method}

\author{Assyr Abdulle\thanks{\email{assyr.abdulle@epfl.ch}} \and Giacomo Rosilho de Souza\thanks{\email{giacomo.rosilhodesouza@epfl.ch}}}

\begin{document}

\maketitle
\blfootnote{\hspace{-4ex}\'Ecole Polytechnique F\'ed\'erale de Lausanne (EPFL), SB-MATH-ANMC, Station 8, 1015 Lausanne, Switzerland.}

\begin{abstract}
An explicit stabilized additive Runge--Kutta scheme is proposed. The method is based on a splitting of the problem in severely stiff and mildly stiff subproblems, which are then independently solved using a Runge--Kutta--Chebyshev scheme. The number of stages is adapted according to the subproblem's stiffness and leads to asynchronous integration needing ghost values. Whenever ghost values are needed, linear interpolation in time between stages is employed. One important application of the scheme is for parabolic partial differential equations discretized on a nonuniform grid. The goal of this paper is to introduce the scheme and prove on a model problem that linear interpolations trigger instabilities into the method. Furthermore, we show that it suffers from an order reduction phenomenon. The theoretical results are confirmed numerically. 
\end{abstract}



\section{Introduction}\label{sec:moddef}
We consider the ordinary differential equation (ODE)
\begin{align}\label{eq:ode}
y'&=f(y) \quad t>0, & y(0)=y_0,
\end{align}
where $y(t)\in\Rn$ with $n\geq 2$, $f:\Rn\rightarrow\Rn$ is a smooth function and $y_0\in\Rn$ is the initial value. We suppose that $f$ can be split in a severely stiff and a mildly stiff term in the following sense: there is a diagonal matrix $D\in\Rb^{n\times n}$ such that $D_{ii}=0$ or $D_{ii}=1$ for $i=1,\ldots,n$ and $Df$, $(I-D)f$ are a severely stiff and a mildly stiff term, respectively.\footnote{We are aware that "severely stiff" and "mildly stiff" are qualitative somewhat imprecise characterizations. This is meant to indicate that the fastest dynamics are in the severely stiff terms. Since the slower scales can still be fast enough to prevent the use of classical explicit schemes, we call them mildly stiff.} A typical example are spatially discretized parabolic problems with a locally refined region. Time discretization leads to a system of ODEs where the eigenvalues of the Jacobian depend on the mesh size and severely stiff components correspond to the refined region. In contrast, mildly stiff components correspond to the coarse region (where the CFL condition still holds).

Since there is a stiff term $Df$ then integration of \eqref{eq:ode} becomes expensive, even if stiffness is induced by a few components only. Multirate methods exploit the special structure of the problem in order to reduce the computational cost. This is often achieved by adapting the Runge--Kutta (RK) method or the step size to the specific partition of the system and employing interpolations or extrapolations for coupling the components together. It is known that the coupling strategy between the stiff and nonstiff terms strongly affects the stability of the system. Indeed a major difficulty in the field is to construct stable multirate methods (see for instance \cite{GeW84,GKR01,GuR93,Hof76,Kva99}). 

The goal of this report is to discuss the properties of an additive Runge--Kutta--Chebyshev (RKC) scheme, which turns out to be very similar to the method described in \cite{Mir17}.
First, we will show that the linear interpolations employed in the scheme might render the integration process unstable. Second, we discuss an order reduction phenomenon observed in numerical experiments.

For these reasons, we introduce in \cite{AGR20} a different multirate RKC scheme, called mRKC, that is free of interpolations, explicit and stable.

The remaining of this report is structured as follows. In \cref{sec:rkcrkc} we introduce the method, show the instabilities on a model problem and discuss the order reduction phenomenon. Numerical experiments are provided in \cref{sec:num}.

\section{The additive Runge--Kutta--Chebyshev method}\label{sec:rkcrkc}
We present here an additive method which uses two RKC schemes simultaneously. Depending on the choice of coefficients the method can be of first- or second-order accurate. The scheme preserves the explicitness of the RKC schemes and does not need any predictor step, but makes use of linear interpolations in time between stages. We will show that these interpolations create instabilities and lead to order reduction.

\subsection{The Runge--Kutta--Chebychev method}\label{sec:rkc}
Chebyshev methods are a family of explicit stabilized Runge--Kutta methods \cite{Abd02,AbM01,GuL60,Leb94,LeM94,Med98,HoS80} with variable number of stages. The number of stages $s$ determines the size of the stability domain, who grows as $\beta s^2$ in the direction of the negative real axis. Methods up to order four have been derived \cite{Abd02}. Among these methods, we consider here the Runge--Kutta--Chebyshev (RKC) methods introduced in \cite{HoS80,Ver80,Ver96,VHS90}. First- and second-order RKC schemes have been derived, for which $\beta\approx 2$ and $\beta\approx 0.65$, respectively.

Let $\tau>$ be the step size, $\rho$ the spectral radius of the Jacobian of $f$ (evaluated in $y_0$) and $s\in\Nb$ such that $\tau\rho\leq \beta s^2$. One step of the RKC scheme is given by
\begin{align}\label{eq:rkc}
\begin{split}
k_0 &= y_0,\\
k_1 &= k_0 +\tau\mu_1 f(k_0),\\
k_i &= \nu_i k_{i-1}+\kappa_i k_{i-2}+(1-\nu_i-\kappa_i)k_0+\tau\mu_i f(k_{i-1})+\tau\gamma_i f(k_0) \quad \mbox{for }i=2,\ldots,s,\\
y_1 &= k_s.
\end{split}
\end{align}
The stages $k_i$ are an approximation of $y(c_i\tau)$, with $\{c_i\}_{i=0}^s$ a strictly increasing sequence satisfying $c_0=0$ and $c_s=1$. The definition of sequences $\mu,\nu,\kappa,\gamma$ and $c$ depends on $s$ and the order of the method.

\begin{figure}[!tbp]
\begin{center}
\begin{subfigure}[t]{\subfigsize\textwidth}
\centering
\includegraphics[trim=2.5cm 0cm 0cm 1cm, clip, width=0.9\textwidth]{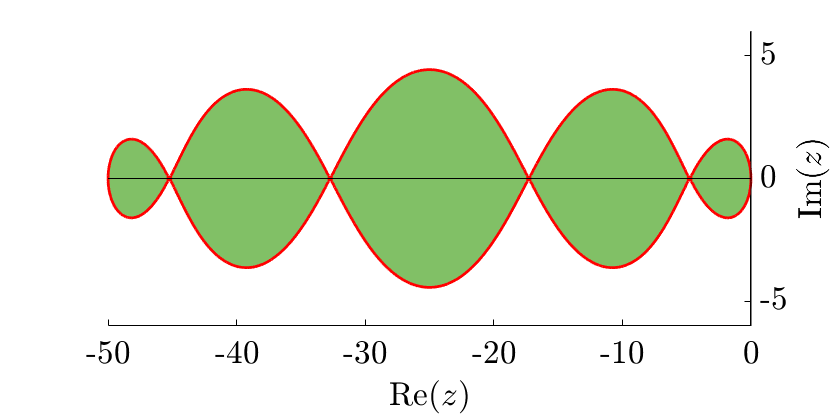}
\caption{Stability domain of the undamped method.}
\label{fig:undamped_stability}
\end{subfigure}
\hfill
\begin{subfigure}[t]{\subfigsize\textwidth}
\centering
\includegraphics[trim=2.5cm 0cm 0cm 1cm, clip, width=0.9\textwidth]{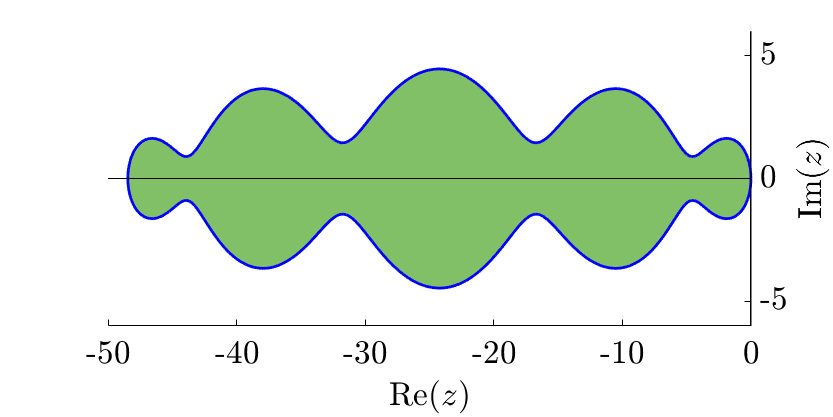}
\caption{Stability domain of the damped method, with $\varepsilon=0.05$.}
\label{fig:damped_stability}
\end{subfigure}
\end{center}
\caption{Stability domain of the damped and undamped first-order RKC method with $s=5$ stages.}
\label{fig:stabrkc1}
\end{figure}

Applying the RKC scheme to the test equation $y'=\lambda y$ with $\lambda\in \mathbb{C}$ one gets $y_1=R_s(\tau\lambda)y_0$, where $R_s$ is the stability polynomial of the RKC scheme. It is clear that $|y_1|\leq |y_0|$ if, and only if, $|R_s(\tau\lambda)|\leq 1$, hence the stability domain of the method is defined by $\mathcal{S}=\{z\in\mathbb{C}\,:\, |R_s(z)|\leq 1\}$. In \cref{fig:undamped_stability} we depict the stability domain of the first-order RKC scheme for $s=5$ and we observe that in some regions the scheme is not stable in the imaginary direction. For this reason, a damping parameter $\varepsilon\geq 0$ is introduced in the method in order to obtain a stability domain containing a narrow strip along the negative real axis (see \cite{GuL60,HoS80} for details). We show in \cref{fig:damped_stability} the stability domain of the first-order RKC scheme with a damping parameter $\varepsilon=0.05$, we observe that it is slightly shorter but stable in the imaginary direction. Taking a damping parameter larger than needed is not convenient since the length of the stability domain decreases and the method would require more function evaluations ($\beta$ is a decreasing function of $\varepsilon$).

\subsection{An additive Runge--Kutta--Chebychev method}
In order to introduce the additive RKC scheme we split \eqref{eq:ode} into a stiff and a mildly stiff problem, which are then integrated independently using two RKC schemes. When communication between the two subproblems is needed, linear interpolation in time is employed.

\subsubsection*{Equation splitting}
Let $D\in \Rb^{n\times n}$ be a diagonal matrix such that $D_{ii}=1$ or $D_{ii}=0$ and $E=I-D$, where $I\in\Rb^{n\times n}$ is the identity matrix. Then \eqref{eq:ode} can be written as 
\begin{align}\label{eq:firstsplit}
y' =& Df(Dy+Ey)+Ef(Dy+Ey)
\end{align}
and multiplying \eqref{eq:firstsplit} either by $D$ or $E$ yields
\begin{subequations}\label{eq:odesplit}
\begin{align}
Dy' =& Df(Dy+Ey),\\
Ey' =& Ef(Dy+Ey),
\end{align}
\end{subequations}
as $DE=ED=0$, $D^2=D$ and $E^2=E$. Letting $y_F=Dy$, $y_S=Ey$, $f_F=Df$ and $f_S=Ef$ we can rewrite \eqref{eq:odesplit} as
\begin{subequations}\label{eq:multirateode}
\begin{align}\label{eq:multirateode1}
y_F' =& f_F(y_F+y_S),\quad t>0 & y_F(0)=Dy_0,\\ \label{eq:multirateode2}
y_S' =& f_S(y_F+y_S),\quad t>0 & y_S(0)=Ey_0.
\end{align}
\end{subequations}
Usually, the matrix $D$ is chosen such that $f_F$ is stiff ($F$ for fast) and $f_S$ is less stiff compared to $f_F$ ($S$ for slow).

\subsubsection*{The additive RKC algorithm}
The additive RKC (ARKC) scheme integrates the two problems in \eqref{eq:multirateode} separately, applying an RKC method to each equation. Integration is performed simultaneously and linear interpolation is employed for the equations' coupling. 

Given $\rhof,\rhos$ the spectral radii of the Jacobians of $\ff,\fs$, respectively, choose $s$ and $m$ such that $\tau\rho_F \leq \beta m^2$ and $\tau\rho_S\leq \beta s^2$. The $\ARKC$ scheme integrates \eqref{eq:multirateode1} using $m$ stages and \eqref{eq:multirateode2} using $s$ stages. Since $\ff$ is supposed to be stiffer than $\fs$ then $m\geq s$.
In the following we call $\mu_i,\nu_i,\kappa_i,\gamma_i$ and $c_i$ the coefficients of an $s$-stage RKC method and $\alpha_j,\beta_j,\delta_j,\zeta_j$ and $d_j$ the coefficients of an $m$-stage RKC method. Further, $k_i$ will be an approximation to $y_S(c_i\tau)$ and $l_j$ an approximation to $y_F(d_j\tau)$. 

If $y_F(t)$ was known, we could integrate \eqref{eq:multirateode2} with the scheme 
\begin{subequations}\label{eq:hypalgo}
\begin{align}\label{eq:hypalgo1}
\begin{split}
k_0 =& y_S(0),\\
k_1 =& k_0 +\tau\mu_1 \fs(y_F(0)+k_0),\\
k_i =& \nu_i k_{i-1}+\kappa_i k_{i-2}+(1-\nu_i-\kappa_i)k_0\\
&+\tau\mu_i f_S(y_F(c_{i-1}\tau)+k_{i-1})+\tau\gamma_i f_S(y_F(0)+k_0) \quad \mbox{for }i=2,\ldots,s.\\
\end{split}
\end{align}
Alternatively, if $y_S(t)$ was known, we could integrate \eqref{eq:multirateode1} with the scheme
\begin{align}\label{eq:hypalgo2}
\begin{split}
l_0 =& y_F(0),\\
l_1 =& l_0 +\tau\alpha_1 \ff(l_0+y_S(0)),\\
l_j =& \beta_j l_{j-1}+\delta_j l_{l-2}+(1-\beta_j-\delta_j)l_0\\
&+\tau\alpha_j \ff(l_{j-1}+y_S(d_{j-1}\tau))+\tau\zeta_j \ff(l_0+y_S(0)) \quad \mbox{for }j=2,\ldots,m.\\
\end{split}
\end{align}
\end{subequations}
However, as neither $y_F$ nor $y_S$ are known they must be approximated. Since 
\begin{align*}
y_F(c_i\tau)\approx &\,y_F(d_{j-1}\tau)+\frac{c_i-d_{j-1}}{d_j-d_{j-1}}(y_F(d_j\tau)-y_F(d_{j-1}\tau)) &\mbox{for}&& d_{j-1}<c_i\leq & d_j
\end{align*}
and $l_j\approx y_F(d_j\tau)$ then we approximate $y_F(c_i\tau)$ by $\tilde l_i$ defined by
\begin{subequations}\label{eq:interpl}
\begin{align}\label{eq:interpl1}
\tilde l_i =& l_{j-1}+\frac{c_i-d_{j-1}}{d_j-d_{j-1}}(l_j-l_{j-1}), &\mbox{where}&& d_{j-1}<c_i\leq & d_j.
\end{align}
A similar strategy is used for $y_S(d_j\tau)$, we approximate it by
\begin{align}\label{eq:interpl2}
\tilde k_j =& k_{i-1}+\frac{d_j-c_{i-1}}{c_i-c_{i-1}}(k_{i}-k_{i-1}), & \mbox{where} && c_{i-1}< d_j\leq & c_i.
\end{align}
\end{subequations}

Replacing in \eqref{eq:hypalgo} the exact values $y_F(c_{i}\tau)$, $y_S(d_j\tau)$ by the approximations $\tilde l_i$, $\tilde \kappa_j$, respectively, yields a fully discrete scheme. Letting $y_0=y(0)$, one step of the $\ARKC$ method is given by
\begin{subequations}\label{eq:algo}
\begin{align}\label{eq:algo1}
\begin{split}
k_0 =& Ey_0,\\
k_1 =& k_0 +\tau\mu_1 \fs(l_0+k_0),\\
k_i =& \nu_i k_{i-1}+\kappa_i k_{i-2}+(1-\nu_i-\kappa_i)k_0\\
&+\tau\mu_i f_S(\tilde l_{i-1}+k_{i-1})+\tau\gamma_i f_S(l_0+k_0) \quad \mbox{for }i=2,\ldots,s\\
\end{split}
\end{align}
and
\begin{align}\label{eq:algo2}
\begin{split}
l_0 =& Dy_0,\\
l_1 =& l_0 +\tau\alpha_1 \ff(l_0+k_0),\\
l_j =& \beta_j l_{j-1}+\delta_j l_{l-2}+(1-\beta_j-\delta_j)l_0\\
&+\tau\alpha_j \ff(l_{j-1}+\tilde k_{j-1})+\tau\zeta_j \ff(l_0+k_0) \quad \mbox{for }j=2,\ldots,m,\\
\end{split}
\end{align}
\end{subequations}
where $\tilde k_j,\tilde l_i$ are defined in \eqref{eq:interpl} and $y_1=k_s+l_m$ is an approximation to $y(\tau)$.

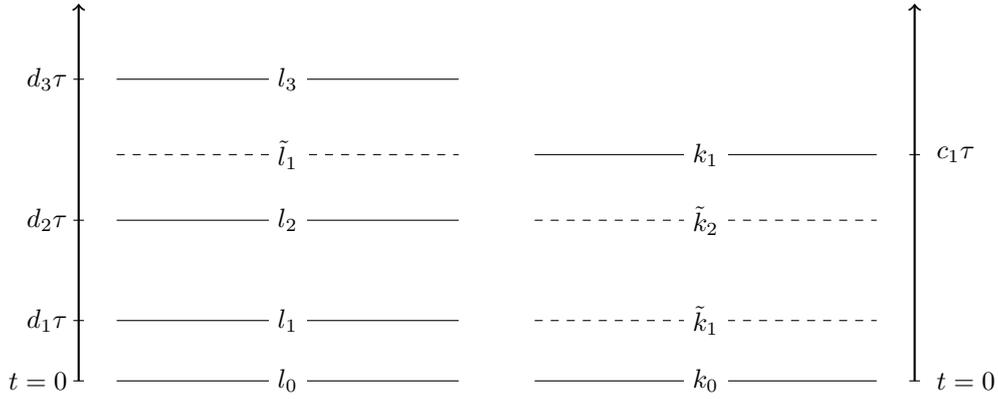
\begin{figure}
\begin{center}
\begin{tikzpicture}
\coordinate (O) at (0,0);
\coordinate (L1) at (5.5,0);
\coordinate (L2) at (-5.5,0);
\coordinate (H1) at (5.5,5);
\coordinate (H2) at (-5.5,5);
\coordinate (d0) at (-5.5,0);   \coordinate (l0a) at (-0.5,0); \coordinate (l0b) at (-5,0);
\coordinate (d1) at (-5.5,0.8); \coordinate (l1a) at (-5,0.8); \coordinate (l1b) at (-0.5,0.8); \coordinate (tk1a) at (0.5,0.8); \coordinate (tk1b) at (5,0.8);
\coordinate (d2) at (-5.5,2.13);\coordinate (l2a) at (-5,2.13);\coordinate (l2b) at (-0.5,2.13);\coordinate (tk2a) at (0.5,2.13);\coordinate (tk2b) at (5,2.13);
\coordinate (d3) at (-5.5,4);   \coordinate (l3a) at (-5,4);   \coordinate (l3b) at (-0.5,4);
\coordinate (c0) at (5.5,0);    \coordinate (k0a) at (0.5,0);  \coordinate (k0b) at (5,0);
\coordinate (c1) at (5.5,3);    \coordinate (k1a) at (0.5,3);  \coordinate (k1b) at (5,3);\coordinate (tl1a) at (-5,3);  \coordinate (tl1b) at (-0.5,3);
\draw (k0a)--(k0b) node [midway,fill=white] {$k_0$};
\draw (k1a)--(k1b) node [midway,fill=white] {$k_1$};
\draw[dashed] (tk1a)--(tk1b) node [midway,fill=white] {$\tilde k_1$};
\draw[dashed] (tk2a)--(tk2b) node [midway,fill=white] {$\tilde k_2$};
\draw[thick,->] (L1)--(H1);% node[anchor=north west] {$c_i\tau$};
\draw (l0a)--(l0b) node [midway,fill=white] {$l_0$};
\draw (l1a)--(l1b) node [midway,fill=white] {$l_1$};
\draw (l2a)--(l2b) node [midway,fill=white] {$l_2$};
\draw (l3a)--(l3b) node [midway,fill=white] {$l_3$};
\draw[dashed] (tl1a)--(tl1b) node [midway,fill=white] {$\tilde l_1$};
\draw[thick,->] (L2)--(H2);% node[anchor=north east] {$d_j\tau$};
\draw ([xshift=-2pt]c0)--([xshift=2pt]c0) node[anchor=west,inner sep=6pt] {$t=0$};
\draw ([xshift=-2pt]c1)--([xshift=2pt]c1) node[anchor=west,inner sep=6pt] {$c_1\tau$};
\draw ([xshift=-2pt]d0)--([xshift=2pt]d0) node[anchor=east,inner sep=6pt] {$t=0$};
\draw ([xshift=-2pt]d1)--([xshift=2pt]d1) node[anchor=east,inner sep=6pt] {$d_1\tau$};
\draw ([xshift=-2pt]d2)--([xshift=2pt]d2) node[anchor=east,inner sep=6pt] {$d_2\tau$};
\draw ([xshift=-2pt]d3)--([xshift=2pt]d3) node[anchor=east,inner sep=6pt] {$d_3\tau$};
\end{tikzpicture}
\end{center}
\caption{Illustration of the ARKC algorithm, solid lines represent the stages $k_i,l_j$ while dashed lines represent the interpolations $\tilde k_j,\tilde l_i$. In this example, the algorithm proceeds as follows. As $k_0$ and $l_0$ are known the scheme can compute $k_1$ and $l_1$, which are approximations at times $c_1\tau$ and $d_1\tau$, respectively, with $d_1<c_1$. Then, it cannot compute $k_2$, as it would need $\tilde l_1$, which is an interpolation of $l_2$ and $l_3$, that are not yet computed. But, it can compute $l_2$ for which $\tilde k_1$, an interpolation between $k_1$ and $k_0$, can be computed. Once $l_2$ is computed, $k_2$ can still not be computed as $l_3$ is missing. But it can compute $l_3$, which requires an interpolation $\tilde k_2$ of $k_1$ and $k_0$ (note that $\tilde k_1,\tilde k_2$ are both interpolations of $k_0$ and $k_1$, but at different times $d_1\tau,d_2\tau$, respectively). Once $l_3$ is known $\tilde l_1$ can be computed and $k_2$ can be evaluated. Informally, we observe that the rule is to advance the variable which is behind in time.}
\label{fig:algo}
\end{figure}

Observe as the conditions on $c_i,d_j$ in interpolations \eqref{eq:interpl} impose an interlaced evaluation order for the stages $k_i$, $l_j$ in \eqref{eq:algo}. For instance, the algorithm can compute both $k_1$, $l_1$ as $k_0$ and $l_0$ are known. But then it can compute $k_2$ only if $c_1\leq d_1$. Indeed, the computation of $k_2$ requires $\tilde l_1$ and the latter needs $l_j$, where $j$ is such that $c_1\leq d_j$. Since only $l_1$ has been computed, the scheme can compute $k_2$ only if $c_1\leq d_1$. Otherwise it computes $l_2$, which can be computed if $d_1\leq c_1$. Hence, at each iteration the algorithm verifies which condition \eqref{eq:interpl1} or \eqref{eq:interpl2} on $c_i,d_j$ is satisfied and computes $k_i$ or $l_j$ accordingly. An illustrative example is provided in \cref{fig:algo}.

The actual implementation of the scheme is fairly simple and a pseudo-code is given in \cref{alg:rkcrkc} below. In the rest of the report we will study how interpolations adversely affect the stability and accuracy of the scheme.
\begin{algorithm}
\caption{$\ARKC$}
\label{alg:rkcrkc}
\begin{algorithmic}
\State Set $s,m$ the smallest integers satisfying $\tau\rho_F\leq \beta m^2$ and $\tau\rho_S\leq \beta s^2$.
\State $k_0=Qy_0$
\State $l_0=Py_0$
\State $k_1=k_0+\tau\mu_1 \fs(l_0+k_0)$ %\Comment First fast step
\State $l_1=l_0+\tau\alpha_1 \ff(l_0+k_0)$ %\Comment First slow step
\State $i=j=1$
\While{$i<s$ or $j<m$} %\Comment Iterate until both $k_i,l_j$ reach the end of the time lapse
\If{$d_j<c_i$} %\Comment Here slow component is late, we compute the next $l_j$
\State $\tilde{k}_j= k_{i-1}+\frac{d_j-c_{i-1}}{c_i-c_{i-1}}(k_{i}-k_{i-1})$ 
\State $j=j+1$
\State $l_{j}=\beta_{j}l_{j-1}+\delta_{j}l_{j-2}+(1-\beta_{j}-\delta_{j})l_0$
\State $\qquad+\tau\alpha_{j} \ff(l_{j-1}+\tilde{k}_{j-1})+\tau\zeta_{j} \ff(l_0+k_0)$
\ElsIf{$c_i\leq d_j$} %\Comment Here fast component is late, we compute the next $k_i$
\State $\tilde{l}_i= l_{j-1}+\frac{c_i-d_{j-1}}{d_j-d_{j-1}}(l_{j}-l_{j-1})$ %\Comment Approx of $y_S(t_0+c_i\tau)$ using slow values at other times
\State $i=i+1$
\State $k_{i}=\nu_{i}k_{i-1}+\kappa_{i}k_{i-2}+(1-\nu_{i}-\kappa_{i})k_0$%\Comment Approx of $y_F(t_0+c_i\tau)$ using RKC
\State $\qquad+\tau\mu_{i} \fs(\tilde{l}_{i-1}+k_{i-1})+\tau\gamma_{i} \fs(l_0+k_0)$
\EndIf
\EndWhile
\State $y_1=k_s+l_m$
\end{algorithmic}
\end{algorithm}

\subsection{Instability}\label{sec:instab}
Now, we study the stability properties of the additive RKC scheme when applied to a $2\times 2$ system. We consider the equation
\begin{align}\label{eq:testeq}
y'=& Ay,\quad t>0 & y(0)=y_0,
\end{align}
where $y_0\in\Rb^2$ and $A\in\Rb^{2\times 2}$ is a symmetric matrix defined by
\begin{align}\label{eq:defA}
A = \begin{pmatrix}
\zeta & \sigma\\ \sigma & \lambda
\end{pmatrix},
\end{align}
with $\lambda,\zeta\leq 0$ and $\sigma^2\leq \lambda\zeta$. Under these conditions $A$ is nonpositive definite. We will study the stability of the $\ARKC$ scheme when applied to \eqref{eq:testeq}. Let
\begin{align*}
D=\begin{pmatrix}
0 & 0\\0 & 1 \end{pmatrix}
\end{align*}
and $E=I-D$. In this setting it holds $f_F(y)=DAy$ and $f_S(y)=EAy$ with $\rho_F=|\lambda|$ and $\rho_S=|\zeta|$. Since the system is linear, applying \cref{alg:rkcrkc} yields
\begin{align*}
y_1 = R_{s,m}(\tau DA,\tau EA)y_0,
\end{align*}
where $s$, $m$ are the number of stages chosen such that $\tau |\lambda|\leq \beta m^2$, $\tau |\zeta|\leq \beta s^2$ and $R_{s,m}(\tau DA,\tau EA)$ is the iteration matrix. The additive RKC scheme is stable if the spectral radius of $R_{s,m}(\tau DA,\tau EA)$ is bounded by one.

Let us fix $s,m\in\Nb$, if $\sigma=0$ then the two equations defined by \eqref{eq:testeq} are independent and the scheme is stable for all $\tau\lambda$ and $\tau\zeta$ such that $\tau|\lambda|\leq \beta m^2$ and $\tau|\zeta|\leq \beta s^2$. We want to investigate the stability of the scheme when $\sigma\neq 0$, hence with coupling. Let $z=\tau\lambda$, $w=\tau\zeta$, $u=\tau\sigma$, then
\begin{align*}
B:=\tau A=\begin{pmatrix}
w & u\\u & z
\end{pmatrix}\qquad\mbox{and}\qquad R_{s,m}(\tau DA,\tau EA)=R_{s,m}(DB,EB).
\end{align*}
Since $\sigma^2\leq \lambda\zeta$ then $u^2\leq zw$ and in the following we consider $u=\theta \sqrt{zw}$ with $\theta\in [-1,1]$. Thus, we define
\begin{align*}
B_\theta = \begin{pmatrix}
w & \theta\sqrt{zw} \\ \theta\sqrt{zw} & z
\end{pmatrix}
\end{align*}
and denote the stability domain of the additive RKC method by
\begin{align*}
\mathcal{S}=\{(z,w)\in \mathbb{R}^2\,:\, \rho(R_{s,m}(DB_\theta,EB_\theta))\leq 1\},
\end{align*}
where $\rho(\cdot)$ denotes the spectral radius of a matrix and $\mathcal{S}$ depends implicitly on the coupling strength $\theta$. We will study the stability of the $\ARKC$ methods for different coupling strengths $\theta$, from $\theta=0$ corresponding to the absence of coupling to $\theta=\pm 1$, the maximal coupling. We will let $(z,w)$ vary in the rectangle $[-\beta m^2,0]\times [-\beta s^2,0]$, which is the stability domain of the method when there is no coupling, i.e. $\theta=0$. The method is considered to be stable if, and only if, $[-\beta m^2,0]\times [-\beta s^2,0]\subset \mathcal{S}$ for all coupling strength $\theta\in [-1,1]$.

\begin{figure}[!tbp]
\begin{center}
\begin{subfigure}[t]{\subfigsize\textwidth}
\centering
\includegraphics[trim=0cm 0cm 0cm 0cm, clip, width=1\textwidth]{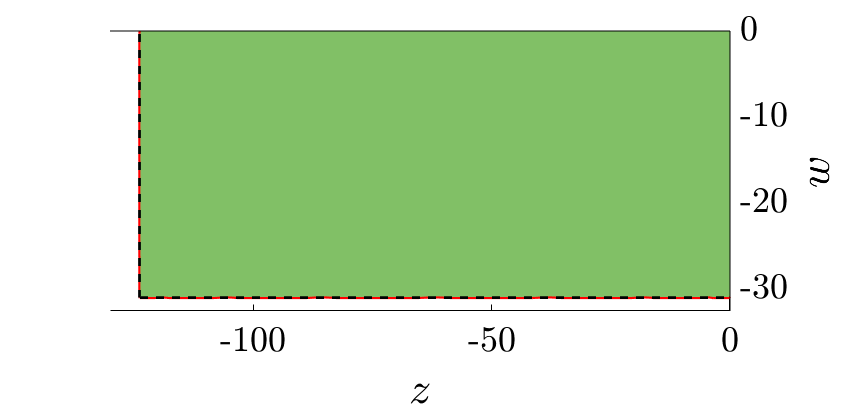}
\caption{Stability domain with $\theta=0$ and $\varepsilon=0.05$.}
\label{fig:stabrkc184_a}
\end{subfigure}
\begin{subfigure}[t]{\subfigsize\textwidth}
\centering
\includegraphics[trim=0cm 0cm 0cm 0cm, clip, width=1\textwidth]{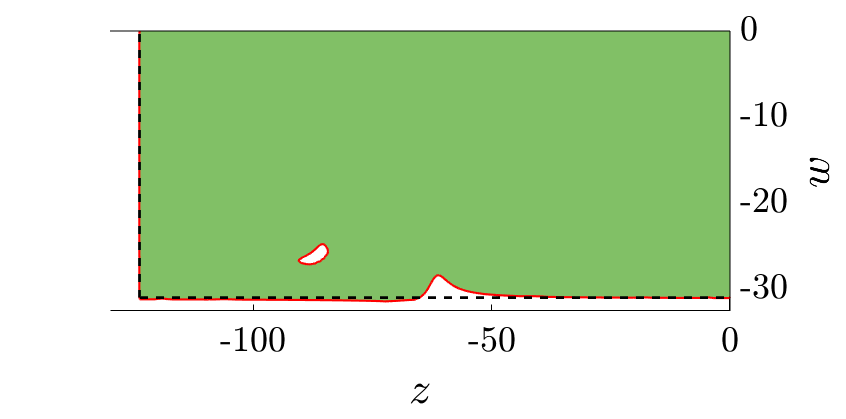}
\caption{Stability domain with $\theta=0.05$ and $\varepsilon=0.05$.}
\label{fig:stabrkc184_b}
\end{subfigure}\\
\begin{subfigure}[t]{\subfigsize\textwidth}
\centering
\includegraphics[trim=0cm 0cm 0cm 0cm, clip, width=1\textwidth]{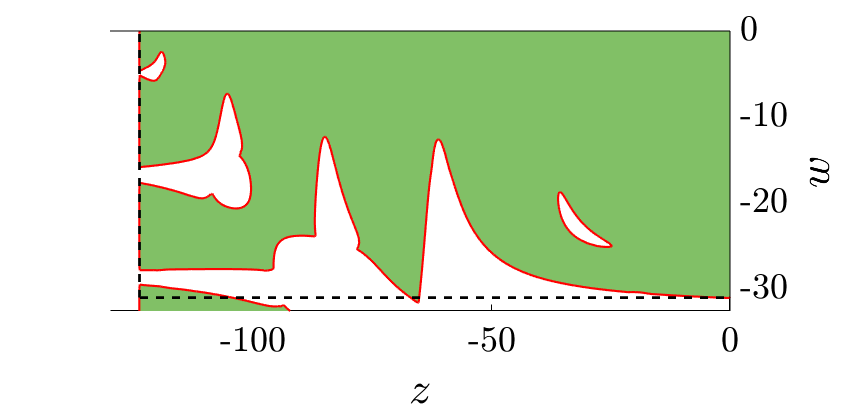}
\caption{Stability domain with $\theta=0.2$ and $\varepsilon=0.05$.}
\label{fig:stabrkc184_c}
\end{subfigure}
\begin{subfigure}[t]{\subfigsize\textwidth}
\centering
\includegraphics[trim=0cm 0cm 0cm 0cm, clip, width=1\textwidth]{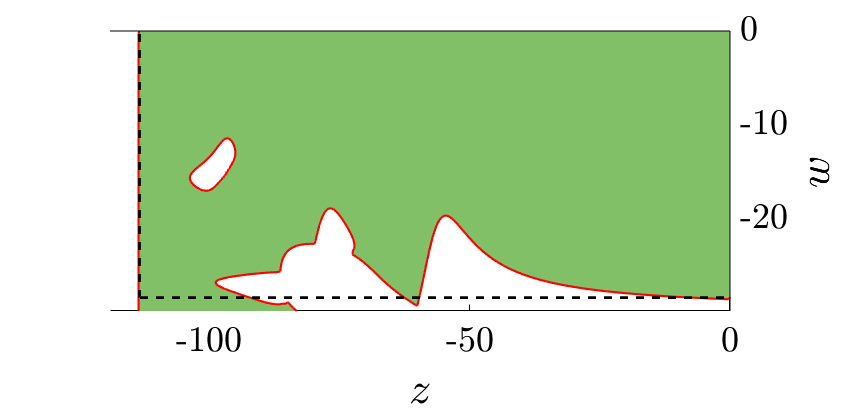}
\caption{Stability domain with $\theta=0.2$ and $\varepsilon=0.2$.}
\label{fig:stabrkc184_d}
\end{subfigure}
\end{center}
\caption{Stability domains of the first-order $\ARKC$ method for $m=8$, $s=4$.}
\label{fig:stabrkc184}
\end{figure}

\begin{figure}[!tbp]
\begin{center}
\begin{subfigure}[t]{\subfigsize\textwidth}
\centering
\includegraphics[trim=0cm 0cm 0cm 0cm, clip, width=1\textwidth]{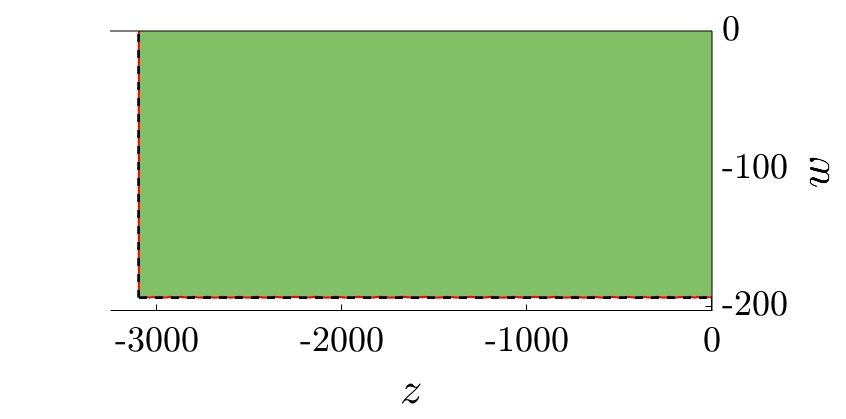}
\caption{Stability domain with $\theta=0$ and $\varepsilon=0.05$.}
\label{fig:stabrkc14010_a}
\end{subfigure}
\begin{subfigure}[t]{\subfigsize\textwidth}
\centering
\includegraphics[trim=0cm 0cm 0cm 0cm, clip, width=1\textwidth]{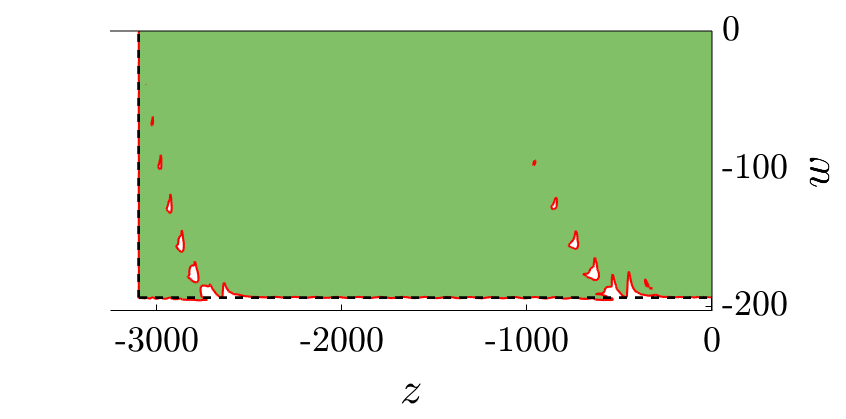}
\caption{Stability domain with $\theta=0.05$ and $\varepsilon=0.05$.}
\label{fig:stabrkc14010_b}
\end{subfigure}\\
\begin{subfigure}[t]{\subfigsize\textwidth}
\centering
\includegraphics[trim=0cm 0cm 0cm 0cm, clip, width=1\textwidth]{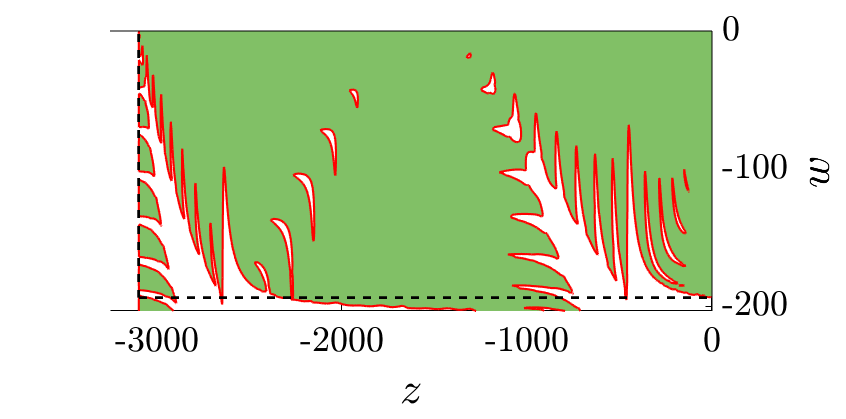}
\caption{Stability domain with $\theta=0.2$ and $\varepsilon=0.05$.}
\label{fig:stabrkc14010_c}
\end{subfigure}
\begin{subfigure}[t]{\subfigsize\textwidth}
\centering
\includegraphics[trim=0cm 0cm 0cm 0cm, clip, width=1\textwidth]{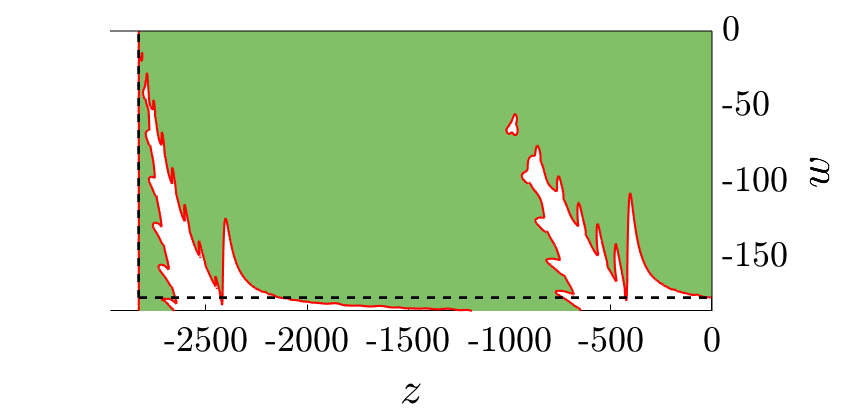}
\caption{Stability domain with $\theta=0.2$ and $\varepsilon=0.2$.}
\label{fig:stabrkc14010_d}
\end{subfigure}
\end{center}
\caption{Stability domains of the first-order $\ARKC$ method for $m=40$, $s=10$.}
\label{fig:stabrkc14010}
\end{figure}

\begin{figure}[!tbp]
\begin{center}
\begin{subfigure}[t]{\subfigsize\textwidth}
\centering
\includegraphics[trim=0cm 0cm 0cm 0cm, clip, width=1\textwidth]{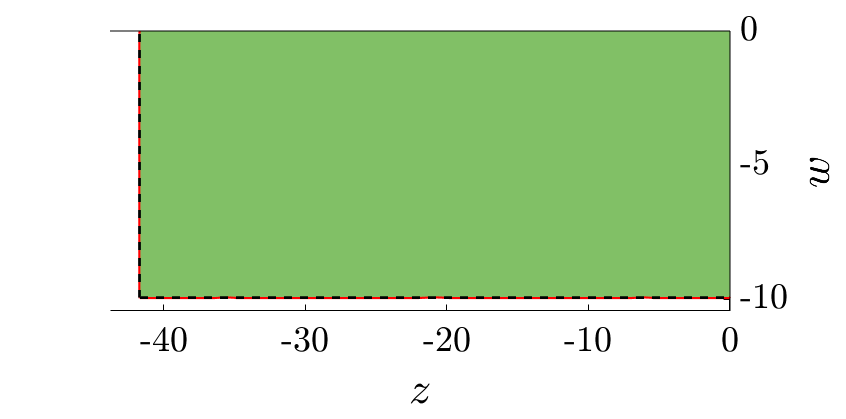}
\caption{Stability domain with $\theta=0$ and $\varepsilon=0.05$.}
\label{fig:stabrkc284_a}
\end{subfigure}
\begin{subfigure}[t]{\subfigsize\textwidth}
\centering
\includegraphics[trim=0cm 0cm 0cm 0cm, clip, width=1\textwidth]{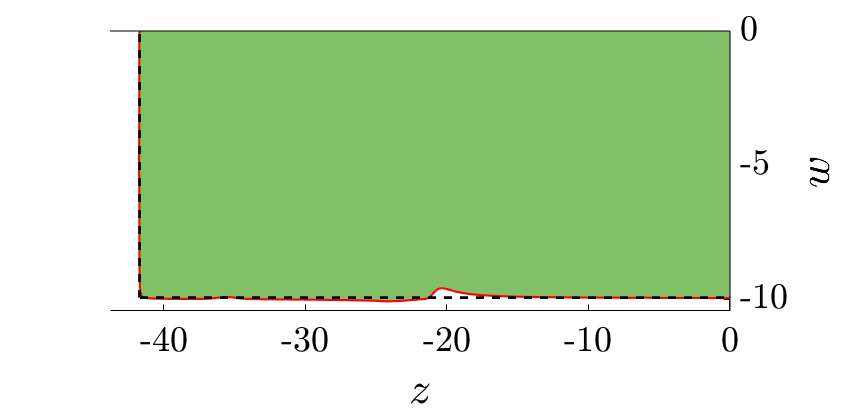}
\caption{Stability domain with $\theta=0.05$ and $\varepsilon=0.05$.}
\label{fig:stabrkc284_b}
\end{subfigure}\\
\begin{subfigure}[t]{\subfigsize\textwidth}
\centering
\includegraphics[trim=0cm 0cm 0cm 0cm, clip, width=1\textwidth]{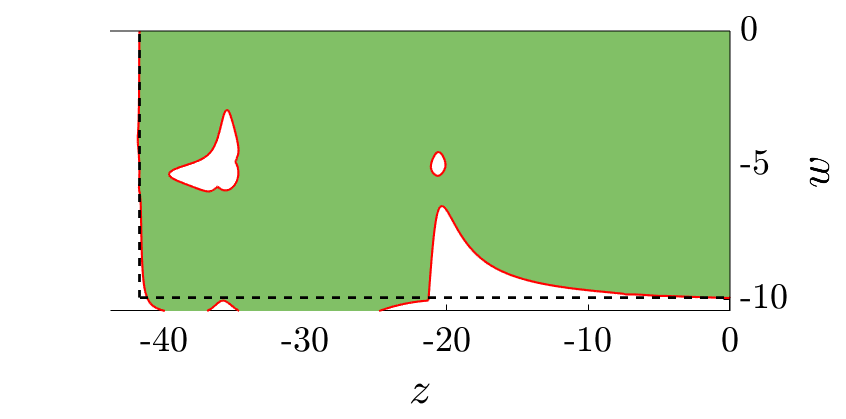}
\caption{Stability domain with $\theta=0.2$ and $\varepsilon=0.05$.}
\label{fig:stabrkc284_c}
\end{subfigure}
\begin{subfigure}[t]{\subfigsize\textwidth}
\centering
\includegraphics[trim=0cm 0cm 0cm 0cm, clip, width=1\textwidth]{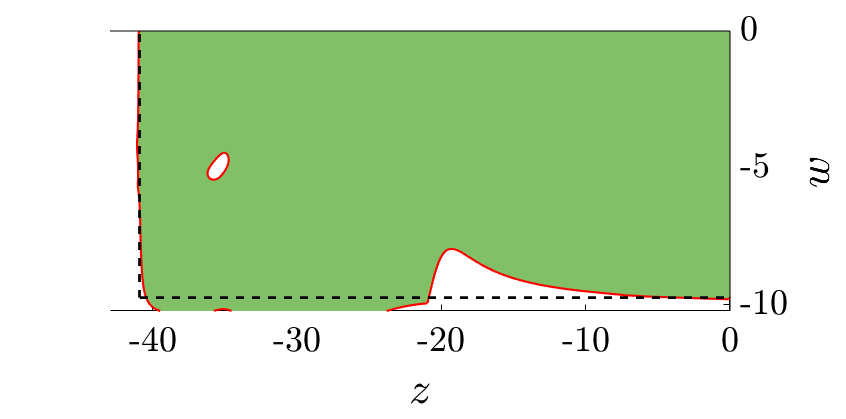}
\caption{Stability domain with $\theta=0.2$ and $\varepsilon=0.2$.}
\label{fig:stabrkc284_d}
\end{subfigure}
\end{center}
\caption{Stability domains of the second-order $\ARKC$ method for $m=8$, $s=4$.}
\label{fig:stabrkc284}
\end{figure}

Observe that the matrix $R_{s,m}(DB_\theta,EB_\theta)$ can be computed replacing $y_0$ by $I$, $f_F(y)$ by $DB_\theta y$ and $f_S(y)$ by $EB_\theta y$ in \cref{alg:rkcrkc}. Hence, for some fixed $s$, $m$ and $\theta$ values we display in \cref{fig:stabrkc184,fig:stabrkc14010,fig:stabrkc284} the stability domain $\mathcal{S}$ of the $\ARKC$ methods by computing the spectral radius of $R_{s,m}(DB_\theta,EB_\theta)$ for varying $(z,w)\in [-\beta m^2,0]\times [-\beta s^2,0]$. The shaded regions represent the stability domains, while the dashed black lines represent the box $B_{s,m}=[-\beta m^2,0]\times [-\beta s^2,0]$, which is the region where the method is stable in absence of coupling, i.e. for $\theta=0$. In \cref{fig:stabrkc184} we show the results for the first-order ARKC method with $m=8$ and $s=4$. We observe in \cref{fig:stabrkc184_a} that for $\theta=0$ and a standard damping parameter $\varepsilon=0.05$, the method is stable in the box $B_{s,m}$, as expected. In \cref{fig:stabrkc184_b,fig:stabrkc184_c} we increase the coupling factor $\theta$ and observe that instability regions appear inside the box $B_{s,m}$. In \cref{fig:stabrkc184_d} we try to increase the damping parameter $\varepsilon=0.2$ and notice that it is not enough to fully stabilize the method. We observed that taking an even larger damping parameter does not stabilize the method. We perform the same experiment in \cref{fig:stabrkc14010} but taking $s=40$ and $m=10$, we see again that if $\theta>0$ the method has instability regions inside the box $B_{s,m}$ and increasing the damping parameter $\varepsilon$ does not help in stabilizing the scheme. Moreover, comparing \cref{fig:stabrkc184,fig:stabrkc14010} we remark that the pattern of the instability region is very different and hence not predictable. We perform the same experiment using the second-order $\ARKC$ method and obtain similar results (see \cref{fig:stabrkc284}). 

\Cref{fig:stabrkc184,fig:stabrkc14010,fig:stabrkc284} illustrate that the additive RKC method discussed here is not stable. Furthermore, the location of the instability regions is not easy to characterize, thus changing the values of $s$ and $m$ does not help in stabilizing the scheme in a given region.

\subsection{Order reduction in the second-order additive RKC scheme}\label{sec:ordred}
For simplicity, we will motivate the order reduction phenomenon using a semi-discrete method, where the stage values $l_j$ are known beforehand and are exact, that is $l_j=y_F(d_j\tau)$. Furthermore, we assume that $\fs$ is integrated exactly by the second-order RKC scheme. In this situation, \cref{alg:rkcrkc} reduces to computing the exact solution $\tilde y_S(\tau)$ of
\begin{align*}
\tilde y_S' =& f_S(\tilde y_F+\tilde y_S),\qquad t\in (0,\tau] & y_S(0)=&Ey_0,
\end{align*}
where $\tilde y_F(t):[0,\tau]\rightarrow \Rb^n$ is the piece-wise interpolation of $y_F(d_j\tau)$ for $j=0,\ldots,m$. 
This scheme is clearly more accurate than the general additive RKC scheme (i.e. when $l_j$ are not known and $\fs$ in not integrated exactly) and an order reduction for this semi-discrete method will imply the order reduction for the full second-order $\ARKC$ method.

%% Showing the standard result
%Let $E_F(t) = y_F(t)-\tilde y_F(t)$ for $t\in [0,\tau]$. By Taylor expansion we can show that for $t=\tilde c\tau \in [c_{i-1}\tau,c_i\tau]$ it holds
%\begin{align*}
%E_F(t)=&y_F(t)-\tilde y_F(t) = y_F(t) -y_F(c_{i-1}\tau)-\frac{\tilde c-c_{i-1}}{c_i-c_{i-1}}(y_F(c_i\tau)-y_F(c_{i-1}\tau)) \\
%=& \frac{1}{2}\frac{(c_i-\tilde c)(c_{i-1}-\tilde c)}{c_i-c_{i-1}}\left(y_F''(\bar t)(c_{i-1}-\tilde c)-y_F''(\tilde t)(c_i-\tilde c)\right)\tau^2,
%\end{align*}
%where $\bar t,\tilde t\in [c_{i-1}\tau,c_i\tau]$. Since $|c_i-\tilde c|\leq |c_i-c_{i-1}|$ and $|c_{i-1}-\tilde c|\leq |c_i-c_{i-1}|$ then
%\begin{align*}
%|E_F(t)| \leq & C_y |c_i-c_{i-1}|^2\tau^2,
%\end{align*}
%where $C_y$ is a constant dependent on $\max_{t\in [0,\tau]} |y''(t)|$.

Let us define $E_F(t) = y_F(t)-\tilde y_F(t)$ and $E_S(t) = y_S(t)-\tilde y_S(t)$ for $t\in [0,\tau]$. For $t \in [d_{j-1}\tau,d_j\tau]$, from a standard linear interpolation result we get
\begin{align*}
|E_F(t)| \leq & C_y ((d_j-d_{j-1})\tau)^2,
\end{align*}
where $C_y$ is a constant dependent on $\max_{t\in [0,\tau]} |y_F''(t)|$. For $t\in [0,\tau]$ it holds
\begin{align*}
E_S(t) =& \int_0^t f_S(y_F(s)+y_S(s))-f_S(\tilde y_F(s)+\tilde y_S(s)) \dif s 
= \int_0^t\int_0^1\frac{\partial f_S}{\partial y}(\bar y(r,s))(E_F(s)+E_S(s))\dif r\dif s,
\end{align*}
with $\bar y(r,s)$ in the segment $[y_F(s)+y_S(s), \tilde y_F(s)+\tilde y_S(s)]$. Supposing $\Vert \frac{\partial f_S}{\partial y}\Vert\leq M_S$ we get
\begin{align*}
| E_S(t)| \leq & M_S\tau \max_{s\in [0,\tau]}| E_F(s)|+M_S\int_0^t | E_S(s)| \dif s
\end{align*}
and using Gronwall's lemma we obtain
\begin{align}\notag
| E_S(t)| \leq & M_S \tau \max_{s\in [0,\tau]}| E_F(s)| e^{\tau M_S} \\ \notag
\leq & C_y M_S e^{\tau M_S} \max_{j=1,\ldots,m}|d_j-d_{j-1}|^2 \tau^3\\ \label{eq:bound0}
=& C_S(\tau) \max_{j=1,\ldots,m}|d_j-d_{j-1}|^2 \tau^3,
\end{align}
where $C_S(\tau)=C_y M_S e^{\tau M_S}$ is bounded from below by $C_y M_S$. 
Let us now estimate the quantity $\max_{j=1,\ldots,m}|d_j-d_{j-1}|^2$ in the nonstiff and the stiff regime.

In a nonstiff regime $\tau\rhof$ is small, where we recall that $\rhof$ is the spectral radius of the Jacobian of $\ff$. Since the stability condition of the second-order RKC method is $\tau\rhof\leq\beta m^2$ then $m$ is a small number in this regime. It follows that the discretization of the interval $[0,1]$ by the nodes $\{d_j\}_{j=1}^m$ is coarse and the estimate
\begin{align}\label{eq:bound1}
\max_{j=1,\ldots,m}|d_j-d_{j-1}|^2\leq 1
\end{align}
is accurate, implying that
\begin{align}\label{eq:nonstiff}
| E_S(t)| \leq C_S(\tau)\tau^3
\end{align}
is tight. Therefore, in a nonstiff regime the interpolation error introduces a third-order local error in the numerical solution, without deteriorating the global second-order accuracy of the $\ARKC$ scheme.

In contrast in a stiff regime $\tau\rhof$ is large and therefore $m$ is large as well. For a damping parameter $\varepsilon=0$ we have $d_j=(j^2-1)/(m^2-1)$ (see \cite{VHS90}) and thus
\begin{align}\label{eq:diffd}
\max_{j=1,\ldots,m}|d_{j}-d_{j-1}|=(d_{m}-d_{m-1})=\frac{m^2-1-(m-1)^2+1}{m^2-1}=\frac{2m-1}{m^2-1} \approx \frac{2}{m},
\end{align}
where the last approximation comes from the fact that $m$ is large. Using $\tau\rhof\leq\beta m^2$ and \eqref{eq:diffd} yield
\begin{align}\label{eq:bound2}
\max_{j=1,\ldots,m}|d_j-d_{j-1}|^2\approx \frac{4}{m^2}\leq \frac{4\beta}{\tau\rhof}.
\end{align}
Hence, from \labelcref{eq:bound0,eq:bound2} we obtain, approximately,
\begin{align}\label{eq:stiff}
|E_S(t)|\leq C_S(\tau)\frac{4\beta}{\rhof}\tau^2.
\end{align}
Let us now discuss both interpolation error \labelcref{eq:nonstiff,eq:stiff}. We observe that in the stiff regime $\tau\rhof\gg 1$, estimate  (2.14) is a second-order interpolation error. In the nonstiff regime $\tau\rhof=\bigo{1}$, estimate \eqref{eq:nonstiff} is accurate and represent a third-order interpolation error. If we now fix $\rhof$ and vary $\tau$ (as done in \cref{fig:conv}) the transition from stiff to nonstiff regime occurs when the step size $\tau$ is sufficiently small so that $\tau\rhof=\bigo{1}$. This is what is seen in \cref{fig:conv}, where a second-order local interpolation error yields a first-order convergence of the $\ARKC$, while for sufficiently small $\tau$ (nonstiff regime) a second-order convergence is recovered thanks to the third-order local interpolation error.

We note that second-order interpolation techniques for the stage values could lead to a genuine third-order interpolation error (also in the stiff regime). However, we observed that second-order interpolations techniques completely destroy the stability of the scheme and we are not aware of a strategy to avoid such instabilities.

\section{Numerical experiments}\label{sec:num}
In this section we present two numerical experiments that support the results of \cref{sec:instab,sec:ordred}.

\subsection{Instability on the model problem}
We show numerically that the first-order $\ARKC$ scheme is unstable (a similar example can be derived for the second-order $\ARKC$ scheme as well).

\begin{figure}
\begin{center}
\begin{tikzpicture}[scale=\plotimscale]
\begin{axis}[height=\aspectratio*\plotimsizeu\textwidth,width=\plotimsizeu\textwidth,legend columns=1,legend style={at={(0.05,0.95)},anchor=north west,draw=none,fill=none},log basis y={2},legend cell align={left},
xlabel={$t$}, ylabel={Solutions $\ell_2$-norm},label style={font=\normalsize},tick label style={font=\small},legend image post style={scale=\legendmarkscale},legend style={font=\small,nodes={scale=\legendfontscale, transform shape}}]
\addplot+[color=NavyBlue,solid,line width=\plotlinewidth pt,mark=x,mark size=\plotmarksizeu pt,mark repeat = 5,mark phase = 0] table [x=t,y=l2N,col sep=comma] 
{data/RKC1-RKC1_sol.csv};\addlegendentry{$\ARKC$, $m=8$ and $s=4$}
\addplot+[color=OrangeRed,line width=\plotlinewidth pt,mark=square,mark size=\plotmarksizeu pt,mark repeat = 5,mark phase = 0] table [x=t,y=l2N,col sep=comma] 
{data/RKC1_sol.csv};\addlegendentry{RKC, $s=8$}
\end{axis}
\end{tikzpicture}
\caption{Stability of first-order additive RKC and standard RKC.}
\label{fig:stab}
\end{center}
\end{figure}
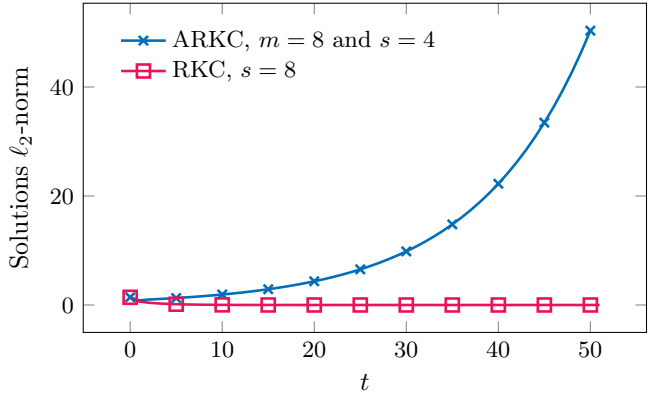

We consider \eqref{eq:testeq} with $y_0=(1,1)^\top$ and $A$ as in \eqref{eq:defA}. We want to choose $s,m,\tau,\lambda,\zeta,\sigma$ such that $s,m$ are the smallest integers satisfying $\tau|\lambda|\leq \beta m^2$, $\tau|\zeta|\leq \beta s^2$ but $\rho(R_{s,m}(\tau PA,\tau QA))>1$. Looking at \cref{fig:stabrkc184_c} we see that the couple $(z,w)=(-100,-28)$ is outside of the stability domain and $m=8$, $s=4$ are the smallest integers satisfying $|z|\leq \beta m^2$ and $|w|\leq \beta s^2$ (recall that $\beta\approx 2$). Hence, if we set $\lambda=-100$, $\zeta=-28$, $\sigma=0.2\sqrt{\lambda\zeta}$, $\tau=1$ and integrate \eqref{eq:testeq} with the $\ARKC$ scheme then it will set $m=8$ and $s=4$. Since $\rho(R_{s,m}(\tau PA,\tau QA))>1$ we expect that the solution explodes. We display in \cref{fig:stab} the norm of the solutions given by the first-order $\ARKC$ and the first-order RKC method, indeed we observe that the $\ARKC$ method is unstable.

\subsection{Order reduction on the heat equation}\label{sec:exp_conv}
\begin{figure}
\begin{center}
\begin{tikzpicture}[scale=\plotimscale]
\begin{loglogaxis}[height=\aspectratio*\plotimsizeu\textwidth,width=\plotimsizeu\textwidth,legend columns=1,legend style={at={(1,0.05)},anchor=south east,draw=none,fill=none},log basis x={2},log basis y={2},legend cell align={left},
xlabel={$\tau$}, ylabel={$\ell_2$-error},label style={font=\normalsize},tick label style={font=\small},legend image post style={scale=\legendmarkscale},legend style={font=\small,nodes={scale=\legendfontscale, transform shape}}]
\addplot+[color=NavyBlue,solid,line width=\plotlinewidth pt,mark=x,mark size=\plotmarksizeu pt] table [x=dt,y=err,col sep=comma] 
{data/order_reduction.csv};\addlegendentry{$\ARKC$}
\addplot[black,domain=0.0625	:0.5] (x,0.5*x);\addlegendentry{$\bigo{\tau}$}
\addplot[black,dotted,thick,domain=0.003953125	:0.0625] (x,2*x^1.5);\addlegendentry{$\bigo{\tau^{1.5}}$}
\addplot[black,dashed,domain=0.00048828125:0.002] (x,x^2/8);\addlegendentry{$\bigo{\tau^2}$}
\end{loglogaxis}
\end{tikzpicture}
\caption{Effective convergence of the second-order additive RKC scheme.}
\label{fig:conv}
\end{center}
\end{figure}
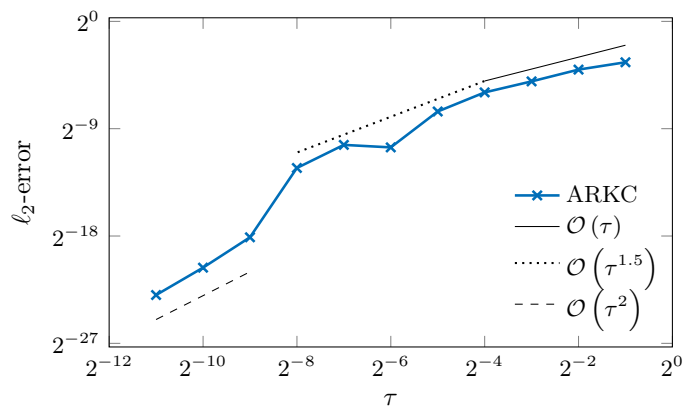

We consider the heat equation
\begin{align}\label{eq:parabolic}
\dt u-\Delta u &= g \qquad\qquad (x,t)\in [0,e]\times [0,1],\\
u(0,t)=u(e,t)&=0 \qquad\qquad\; t\in [0,1],\\
u(x,0)&= u_0(x)\quad\quad\; x\in [0,e],
\end{align}
with $g$, $u_0$ such that the exact solution is $u(x,t)=e^{-t}x(\log(x)-1)$. We discretize the domain $\Omega=[0,e]$ with second-order finite differences. Since $u$ has a spatial singularity in $x=0$ we use a uniform mesh size $H\approx 1/2^4$ in $\Omega_S=(0.005e,e)$ and a uniform mesh size $h\approx H/200$ in $\Omega_F=(0,0.005e)$. After discretization, \eqref{eq:parabolic} can be written as
\begin{align}\label{eq:sdpar}
y' =& Ay+F\quad t\in (0,1], & y(0)=y_0.
\end{align}
Let $D$ be a diagonal matrix of the same size as $A$ such that $D_{ii}=1$ if the $i$th node is in $\overline\Omega_F$ and $D_{ii}=0$ else. We define $f_F(t,y)=D(Ay+F(t))$ and $f_S(t,y)=(I-D)(Ay+F(t))$. We verify the effective order of convergence of the second-order $\ARKC$ scheme integrating \eqref{eq:sdpar} using different step sizes $\tau=1/2^k$, with $k=1,\ldots,11$, comparing the numerical solution against a reference solution computed on the same mesh. We do not use the exact solution since we are only interested in time discretization errors. The results are shown in \cref{fig:conv}, we observe that for $\tau$ large enough the rate of convergence is one, then there is a transition phase and finally for $\tau$ very small the second-order convergence rate is recovered. This result is in line with the findings of \cref{sec:ordred}.

\section{Conclusion}
In this report we discussed an additive Runge--Kutta--Chebyshev method for multirate ordinary differential equations. The method is based on a decomposition of the original problem in two subproblems and integrates both problems with a Runge--Kutta--Chebyshev method, where the number of stages is adapted to the stiffness (fastest scale) of each subproblem. The different stages number leads to an asynchronous integration procedure and linear interpolation in time between stages is employed whenever coupling values are needed.

The scheme is explicit and straightforward to implement. However, we have shown on a model problem that linear interpolations might render the scheme unstable. Furthermore, the second-order additive Runge--Kutta--Chebyshev method suffers from an order reduction phenomenon. Numerical examples corroborate the theoretical findings.

\section*{Acknowledgements} The authors are partially supported by the Swiss National Science Foundation, under grant No. $200020\_172710$.

%\bibliographystyle{abbrv}
%\bibliography{../../../../LaTeX/library}

\end{document}